\documentclass[11pt]{article}
\title{Study on the family of $K3$ surfaces induced from the
lattice $(D_4)^3\oplus \langle -2\rangle \oplus \langle 2\rangle$}
\author{K.Koike, H.Shiga, N.Takayama and T.Tsutsui}
\date{Mar. 18, 2000}

\usepackage{amsfonts}
\usepackage{amsthm}
\usepackage{amssymb}
\usepackage{amsmath}
\usepackage{amscd}
\usepackage[dvips]{graphicx}

\makeatletter
  
      \@addtoreset{equation}{section}
\makeatother
\begin{document}
\maketitle

\def\bigzerol{\smash{\hbox{\b 0}}}
\def\bigzerou{\smash{\lower1.7ex\hbox{\b 0}}}
\newtheorem{thm}{Theorem}
\newtheorem{cor}{Corollary}[section]
\newtheorem{df}{Definition}[section]
\newtheorem{lem}{Lemma}[section]
\newtheorem{rem}{Remark}[section]
\newtheorem{prop}{Proposition}[section]


\setcounter{section}{-1}
\section{Introduction}
\quad \ In 1992 K.Matsumoto, T.Sasaki and M. Yoshida \cite{MSY} studied the 
period mapping
for a family of $K3$ surfaces of type $(3,6)$, that is the family of 
double sextic surface over ${\bf P}^2$ ramified along 6 lines in general 
position, and Matsumoto \cite{Ma} gave the description
of the inverse mapping in terms of theta constants. 
It gives the modular mapping for the 4 dimensional Shimura variety
in the Siegel upper half space ${\cal H}_4$
derived from the family of 4-dimensional abelian varieties with generalized
complex multiplication by $\sqrt{-1}$. So we call it MSY modular mapping.
Shiga showed an arithmetic application of MSY modular mapping in \cite{Sh}. 
This story is the consequence of
the eventual coincidence of 2 different bounded symmetric domains between
$D_{IV}^4$ and $I_{2,2}$. There are a few (finite) such exceptional coincidences.
The highest one is the (analytic) equivalence between $D_{IV}^6$ and $D_{II}(4)$
( in terms of Lie algebra $so(2,6;{\bf R})\cong so(4,{\bf H})$)
 and it contains the above coincidence of MSY case.
\par \indent
Our present study is the first step to get the extended model of the
MSY modular mapping using this equivalence.
\par
Let us consider the rank 14 lattice 
$P=D_4^3\oplus \langle -2\rangle \oplus \langle 2\rangle$.
We define a $K3$ surface $S$ of type $P$ with 
the property that
$P\subset {\rm Pic}(S) $ (see Definition 2.1),
where ${\rm Pic}(S) $ indicates the Picard lattice of $S$.
In this article we study the family of $K3$ surfaces of type $P$ with a 
certain fixed multi-polarization. 
And we do not discuss the representation of the inverse of the period
mapping.
\par
We mention that our family is already appeared in the work \cite{MS}
( in section 7) and
they obtained the differential equation coming from this family standing
on a different view point.
\par \vskip 3mm
\indent
Throughout this article we work on the field ${\bf C}$.
\par
In Section 1 we study the family ${\cal F}$ of double covering surfaces
over ${\bf P}^1\times {\bf P}^1$ branching along 4 bidegree $(1,1)$
curves imposed with a certain generality condition. 
The element of ${\cal F}$ is called a double $4H$ surface ( see Definition 2.2).
Such a surface is given in the affine form
\[
S=S(x): w^2=\prod_{k=1}^4
(x_{1}^{(k)}st+x_{2}^{(k)}s+x_{3}^{(k)}t+x_{4}^{(k)}),
\]
\par
where we use the notation
\[
x_k=\begin{pmatrix}x_{1}^{(k)}&x_{2}^{(k)}\cr x_{3}^{(k)}&x_{4}^{(k)}
\end{pmatrix}
\in M(2,{\bf C}).
\]
We set up the view point that a general member $S$ of ${\cal F}$ 
to be an elliptic
fibred surface (Proposition 1.1).
It becomes to be a $K3$ surface.  Then we construct the basis
of the transcendental lattice ${\rm Tra}(S)$ of a general member $S$.
Fortiori we know
the structures of the Picard lattice ${\rm Pic}(S)$: 
\par
\subsection*{Theorem 1} \quad For a general member $S$ of ${\cal F}$ it holds
\[
{\rm Tra}(S)\cong U(2)^2\oplus \langle -2\rangle ^4.
\]
By making its orthogonal complement in the $K3$ lattice $L$ we get
the Picard lattice 
\[
{\rm Pic}(S)\cong D_4^3\oplus \langle -2\rangle \oplus \langle 2\rangle .
\]
\par
In Section 2 we show that a $K3$ surface $X$ of type $P$ 
can be realized as a double covering surface studied in Section 1 provided
a certain fixed polarization (Theorem 2).
\par
In Section 3 we define a fixed marking of $S$ of a $K3$ surface of type $P$.
First we study the period domain for the family of such marked surfaces.  
That is a 6 dimensional domain given in the form
\[
D^+=\{ \eta =[\eta _1,\cdots , \eta _8]\in {\bf P}^7:
{}^t\eta A\eta =0,{}^t{\overline \eta }A\eta >0, \Im (\eta _3/\eta _1)>0\} , 
\]
where $A=U(2)\oplus U(2)\oplus (-2I_4)$. It is a bounded symmetric domain of type IV.
Next we determine the modular group for the isomorphism classes of the
marked surfaces. It is given as
the principal congruence subgroup $G^+(2)$ of level 2 in the 
full group 
\[
G^+=\{ g\in PGL(8,{\bf Z}): {}^tgAg=A, g(D^+)=D^+\}
\] 
of the positive isometries for the lattice $U(2)^2\oplus \langle -2\rangle ^4$. 
The exact statement is given in Theorem 3.
\par
In Section 4
we consider a general form of the period
\[
u(x)=\int _{C}\Omega =\int _C\prod_{k=1}^4
(x_{1}^{(k)}st+x_{2}^{(k)}s+x_{3}^{(k)}t+x_{4}^{(k)})^{-1/2}ds\wedge dt,
\]
where $C$ indicates an element of $H_2(S(x),{\bf Z})$. We describe
the differential equation for $\int _{C}\Omega $ as a 
multi-valued analytic function of 
$16$ variables $x=(x_1,\cdots ,x_4)$, where we use the notation
\[
x_k=\begin{pmatrix}x_{1}^{(k)}&x_{2}^{(k)}\cr x_{3}^{(k)}&x_{4}^{(k)}\end{pmatrix} \in M(2,{\bf C}).
\]
That is a certain type of GKZ hypergeometric differential equation ( Proposition 4.1 and 4.2) that 
is not so called Aomoto-Gelfond type ( compare with \cite{MS}). We show
the regularity of the system outside some divisor in the
parameter space as stated in Theorem 4.  
We show the regular holonomicity of the system  and 
make the calculation of the rank
based on the theory of \cite{SST-BOOK} .
The exact statement is given in Theorem 5.
We determine the monodromy group for our  system in Theorem 6.
\par
In Section 5 we construct the family $\{ KS(\eta ):\ \eta \in D^+\} $ of
 Kuga-Satake varieties 
corresponding to 
${\cal F}$. That becomes to be equivalent with the family of $8$ 
dimensional Abelian varieties
with the Hamilton quaternion endomorphism structure
(Theorem 7).

\section{Lattice structure of a double $4H$ surface}
\subsection{Setting up the situation}
\quad We consider an algebraic surface $S'$ obtained as a double cover over
${\bf P}^1\times {\bf P}^1$ ramifying along 4 different rational curves 
$H_1,H_2,H_3,H_4$ of bidegree $(1,1)$. Here we suppose the following 
generality condition:
\par
(g1)\quad $H_i\ \ (i=1,2,3,4)$ is irreducible,
\par
(g2)\quad $H_i\cap H_j\ \ (i\not= j)$ consists of 2 different points,
\par
(g3)\quad For any different 3 indices $i,j,k$ we have 
$H_i\cap H_j\cap H_k=\emptyset $.
\\ \vskip 1mm
We denote $\pi $ the projection $S'\rightarrow {\bf P}^1\times {\bf P}^1$.
Set $L_i=\pi^{\ast } H_i$.
The surface $S'$ has 12 singular points of type $A_1$ corresponding to 
the intersections 
$L_i\cap L_j\enskip (i\not= j)$. By the desingularization procedure we get
a $K3$ surface $S$. \par
If we have the condition \par
(e1) \quad the algebraic variety $S'$ has at most simple singularities
\\
instead of $\mathrm{(g1),(g2),(g3)}$, we obtain a $K3$ surface by the
same procedure.
Henceforth we describe the curve $H_k = H(x^k) \ (k =1,2,3,4)$ in the form
\[ (s, 1) \begin{pmatrix}x_1^{(k)} & x_2^{(k)} 
   \\ x_3^{(k)} & x_4^{(k)}
         \end{pmatrix} \begin{pmatrix}t\\1\end{pmatrix} = 0\]
with a matrix
\[ x^k = \begin{pmatrix}x_1^{(k)} &x_2^{(k)} \\ x_3^{(k)} &x_4^{(k)}
         \end{pmatrix} \in \mathrm{M}(2,{\bf C})\]
and an affine coordinate $(s,t)$ of ${\bf P}^1\times {\bf P}^1$.
Let $S'=S'(x)$ be the algebraic variety obtained from 
$ H(x^1)\cup \cdots \cup H(x^4)$ by the above way, and let $S(x)$ be
its desingularization.
Set 
\[ X^{\circ} = \{ x = (x^1,x^2,x^3,x^4) \in \mathrm{M}(2,{\bf C})^4 :
H(x^1),\cdots , H(x^4) \  \text{satisfy} \  (g1),(g2),(g3) \},\]
and set
\[ X' = \{ x = (x^1,x^2,x^3,x^4) \in \mathrm{M}(2,{\bf C})^4:
S' \  \text{has at most simple singularities}\}.\]

\begin{df}
We call {\em a double $4H$ surface}
the $K3$ surface $S$ obtained as $S=S(x), x \in X^{\circ}$.  
{\em An extended double $4H$ surface} is a $K3$ surface obtained as 
$S=S(x), x \in X'$.
Let ${\cal F}$ denote the totality of double $4H$ surfaces:
\[ \{S(x) : x \in X^{\circ} \}\]
We use the following notations.
\par
L: the $K3$ lattice $E_8^2\oplus U^3$,where $E_8$ denotes the negative 
definite even unimodular lattice of rank 8 and $U$ denotes the
 hyperbolic lattice $\begin{pmatrix}0&1\\1&0\end{pmatrix}$. The basis of 
 $L$ is fixed once and always.
\par
${\rm Pic}(S)$: the Picard lattice of $S$, 
\par
${\rm Tra}(S)$: the transcendental lattice of $S$. That is defined as
the orthogonal complement of ${\rm Pic}(S)$ in $H^2(S,{\bf Z})$  
\end{df}
\begin{rem}
We consider $X^{\circ}$ to be the parameter space of $\mathcal{F}$, but
 it contains many abundant parameters. Our family $\mathcal{F}$ contains 
 6 essential parameters. We shall discuss the problem of abundance in
 Section 2 and Section 3.
\end{rem}
\subsection{Elliptic fibration}
\quad In this section we determine ${\rm Tra}(S)$ and ${\rm Pic}(S)$ for a 
general member of ${\cal F}$. For this work we always consider the problem
in the dual lattice $H_2(S,{\bf Z})$ and we observe a special member
of ${\cal F}$.
\par
Let $s,t$ be the affine coordinates of ${\bf P}^1\times {\bf P}^1$.
Henceforth we denote the $s$-space ($t$-space ) by 
${\bf P}^1(s) ({\bf P}^1(t))$, respectively.
Set
\[
f_1(s)=s, \quad f_2(s)=\frac{16}{s}, \quad f_3(s)=\frac{5s+2}{2s+5}, 
\quad f_4(s)=\frac{3s+64}{4s+27}
\]
Consider the following four $(1,1)$ curves :
\[
H_i: t=f_i(s)\enskip (i=1,2,3,4),
\]
and make the double $4H$ surface : 
\[
S_0: w^2=\prod_{i=1}^4(t-f_i(s))
\]
that is derived from the above system
 $\{ H_1,\cdots ,H_4\} $.
Let $\pi_1 $ be the projection from $S_0$ to ${\bf P}^1(s)$.
For a generic point on ${\bf P}^1(s)$ we obtain an elliptic curve $\pi ^{-1}(s)$.
So we get an elliptic surface $(S_0,\pi_1 ,{\bf P}^1)$.
\par
\begin{figure}[hbtp]
 \begin{center}
 \includegraphics[width=6cm]{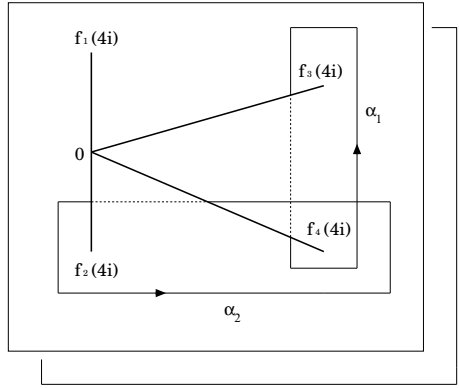}
 \caption{double covering of ${\bf P}^1(t)$}
 \label{fig:dc}
 \end{center}
\end{figure}

 It is easy to observe the following.
\begin{prop}
The elliptic surface $(S_0,\pi _1,{\bf P}^1)$ has 12 singular 
fibres corresponding to the intersections 
$H_i\cap H_j\enskip (i\not= j)$. Every singular fibre is of type $I_2$ 
( according to the Kodaira classification of the singular fibres).
These are situated over real
points :
\begin{eqnarray*}
&&
s=\pm 4 \quad (=H_1\cap H_2),
\pm 1 \quad (=H_1\cap H_3), 
2 \ {\rm and}\ -8 \quad (=H_1\cap H_4),
\cr && 
-2 \ {\rm and}\ 8 \quad (=H_2\cap H_3), 
\pm 12 \quad (=H_3\cap H_4), 
\pm \sqrt{19} \quad (=H_3\cap H_4).
\end{eqnarray*} 
We denote these 12 points by $s_j\ \ (j=1,\cdots ,12)$ according to the
ascending order and set ${\rm Sing}=\{ s_1,\cdots ,s_{12}\} $.
\end{prop} 
We fix a base point $s_0=4\sqrt{-1}$ and set $E_0=\pi_1 ^{-1}(s_0)$. 
It is given by 
\[
w^2=(t-f_1(4i))(t-f_2(4i))(t-f_3(4i))(t-f_4(4i)).
\]
So this is a double cover over ${\bf P}^1(t)$ branched at 4 points
$t=f_k(4i)\ (k=1,2,3,4)$. On $E_0$ we take
a basis $\alpha _1,\alpha _2$ of $H_1(E_0,{\bf Z})$ 
with the intersection multiplicity
$\alpha _1\cdot \alpha _2 =1$ as in Fig.\ref{fig:dc}.
\par
\begin{figure}[hbtp] 
 \begin{center}
 \includegraphics[width=6cm]{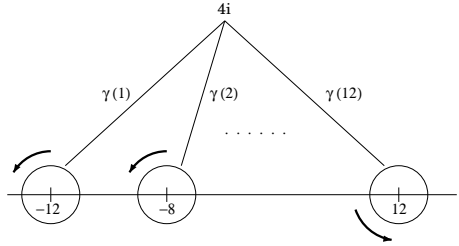}
 \end{center}
 \caption{generator of $\pi_1({\bf P}^1(s)-Sing)$}
 \label{fig:gen}
\end{figure}

The fundamental group $\pi _1({\bf P}^1-{\rm Sing},\ast )$ acts on
$H_1(E_0,{\bf Z})$ as the monodromy group. We describe the generator system
for it.
Let $\gamma _i\enskip (i=1,\cdots ,12)$ be a loop
starting from $s_0$ and goes around $s_i$ in the positive sense 
on the upper half plain except the circuit around $s_i$ 
(see Fig.\ref{fig:gen} ).
By a direct observation we can calculate the monodromy transformation 
$T_j=T(\gamma _j)$
of $H_1(E_0,{\bf Z})$ along $\gamma _j$
with respect to the basis $\{ \alpha _1,\alpha _2\} $.
\begin{prop}
The circuit matrices for $T_j=T(\gamma _j)$ are given by the following table. 
Here $T_j$ acts from right to the system $\{ \alpha _1,\alpha _2\} $. \\
\begin{center}
 \begin{tabular}{|c||c|c|c|c|}
 \hline
 $s_j$ & $\pm4, \pm\sqrt{19}$ & $\pm1, \pm12$ & $8, -2$ & $2, -8$ \\ 
 \hline
 $T_{j}$ & $\begin{pmatrix} 1&2\\ 0&1 \end{pmatrix}$
        & $\begin{pmatrix} 1&0\\ -2&1 \end{pmatrix}$ 
        & $\begin{pmatrix} -1&2\\ -2&3 \end{pmatrix}$
        & $\begin{pmatrix} 3&2\\ -2&-1 \end{pmatrix}$ \\
 \hline
 vanishing cycle & $\alpha_1$ & $\alpha_2$ & $\alpha_1 + \alpha_2$ 
  & $\alpha_1- \alpha_2$ \\
 \hline
 \end{tabular}\end{center}
\end{prop}
\par \vskip 3mm
We choose a fixed point ${\underline s}$ in the lower half plane of
${\bf P}^1(s)$. Make cut lines $c_i$ by the line segment between ${\underline s}$ and 
$s_i\ \ (i=1,\cdots ,12)$. By restricting $(S_0,\pi _1,{\bf P}^1)$ on
${\bf P}^1-\cup _{i=1}^{12}c_i$ we have a topologically trivial fibration.
So we can determine the cycles $\alpha _1,\alpha _2$ of 
$H_1(\pi ^{-1}(s),{\bf Z})$ for
any $s\in {\bf P}^1-\cup _{i=1}^{12}c_i$ using this trivialization.
If we make the continuation of the system $\{ \alpha _1,\alpha _2\}$ 
passing through the line $c_i$ from left to right,
It is transformed according to Table 1.
\subsection{Two systems $\{ \Gamma _i\} $ and $\{ C_i\} $ of 
$H_2(S_0,{\bf Z})$}
\begin{figure}[hbtp]
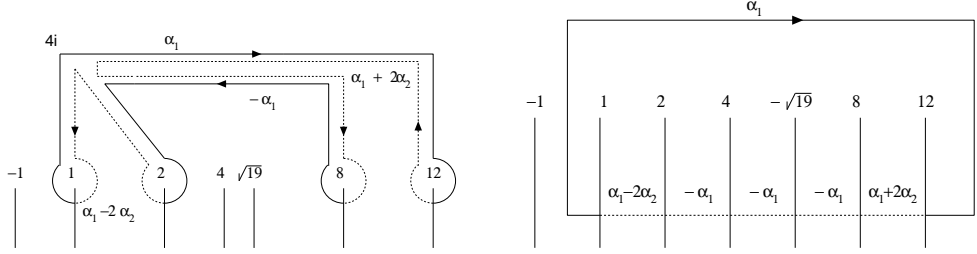

 \begin{center} 
  \includegraphics[width=6cm]{GG1.eps} \qquad
  \includegraphics[width=6cm]{G1.eps}
\end{center}
 \caption{$\Gamma_1$(these are homologically equivalent)}
 \label{fig:ga1}
\end{figure}

\begin{figure}[hbtp]
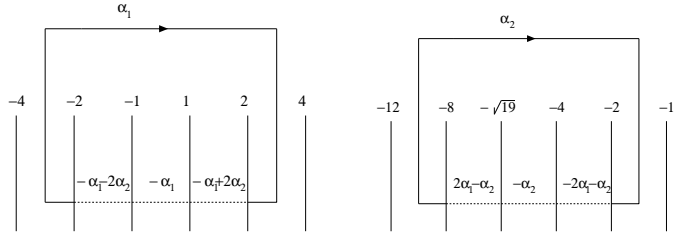

 \begin{center} 
 \includegraphics[width=4cm]{G3.eps} \qquad
 \includegraphics[width=4cm]{G4.eps}
 \end{center}
 \caption{$\Gamma_3$(left),\quad $\Gamma_4$(right)}
 \label{fig:ga34}
\end{figure}

\begin{figure}[hbtp]
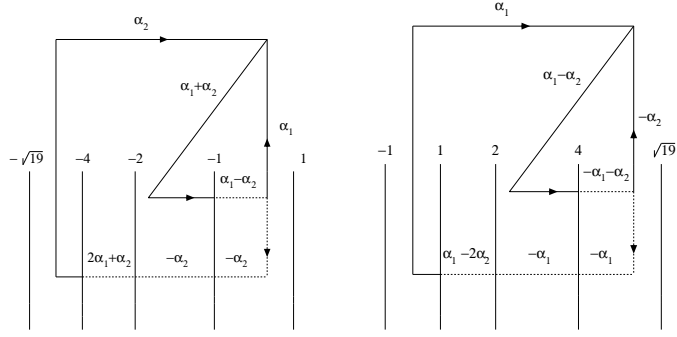

 \begin{center} 
 \includegraphics[width=4cm]{G6.eps} \qquad
 \includegraphics[width=4cm]{G7.eps}
 \end{center}
 \caption{$\Gamma_6$(left), \quad $\Gamma_7$(right)}
 \label{fig:ga67}
\end{figure}

\quad Let $r$ be an oriented arc on ${\bf P}^1(s)$ starting from $s_0$. 
We make a 2-chain on $S_0$ by the continuation of a starting 1-cycle 
$\alpha \in H_1(E_0,{\bf Z})$ along $r$, and denote it by
\[ r\times \alpha . \]
Here we define its orientation as the ordered pair of the ones of $r$ and 
$\alpha $.
Using this notation we make a system
$\{ \Gamma _1,\cdots , \Gamma _8\} $ of 8 elements in $H_2(S_0,{\bf Z})$
(see Fig.\ref{fig:ga1}, Fig.\ref{fig:ga2}, Fig.\ref{fig:ga34}, 
 Fig.\ref{fig:ga67}):
\begin{eqnarray*}
&&
\Gamma _1=\gamma(7)^{-1}\gamma(8)^{-1}\gamma(9)^{-1}\gamma(11)^{-1}\gamma(12)^{-1}\times \alpha _1,
\cr &&
\Gamma _2=\gamma(4)^{-1}\gamma(5)^{-1}\gamma(6)^{-1}\gamma(7)^{-1}\gamma(8)^{-1}\gamma(9)^{-1} 
\times \alpha _2,
\cr &&
\Gamma _3=\gamma(5)^{-1}\gamma(6)^{-1}\gamma(7)^{-1}\gamma(8)^{-1}\times \alpha _1,
\cr &&
\Gamma _4=\gamma(2)^{-1}\gamma(3)^{-1}\gamma(4)^{-1}\gamma(5)^{-1}\times \alpha _2,
\cr &&
\Gamma _5=\gamma(12)\times (\alpha _1+\alpha _2)+\gamma(10)^{-1}\gamma(11)^{-1}\gamma(12)^{-1}\times (-\alpha _2),
\cr &&
\Gamma _6=\gamma(6)\times (\alpha _1+\alpha _2)+\gamma(4)^{-1}\gamma(5)^{-1}\gamma(6)^{-1}\times (-\alpha _2),
\cr &&
\Gamma _7=\gamma(9)\times (\alpha _1-\alpha _2)+\gamma(7)^{-1}\gamma(8)^{-1}\gamma(9)^{-1}\times (-\alpha _1),
\cr &&
\Gamma _8=\gamma(3)\times (\alpha _1-\alpha _2)+\gamma(1)^{-1}\gamma(2)^{-1}\gamma(3)^{-1}\times (-\alpha _1).
\end{eqnarray*}
where the composite arcs are performed from right to left.

\begin{figure}[hbtp]
 \begin{center} 
 \includegraphics[width=6cm]{G2.eps}
  \end{center}
 \caption{$\Gamma_2$}
 \label{fig:ga2}
\end{figure}

We note that any of $\{ \Gamma _1,\cdots , \Gamma _8\} $
is a 2-cycle according to the monodromy action in Table 1. So we
regard them a system in $H_2(S_0,{\bf Z})$.
Let $L(\Gamma )$ denote the sublattice
$\oplus {\bf Z}\Gamma _i$.
\par
\begin{figure}[hbtp]
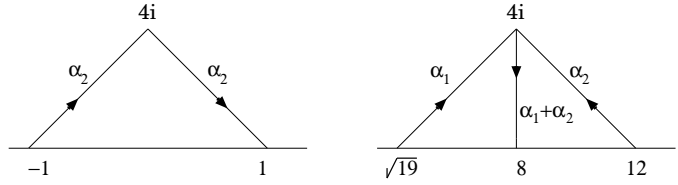

 \begin{center} 
 \includegraphics[width=4cm]{C1.eps} \qquad
 \includegraphics[width=4cm]{C5.eps}
 \end{center}
 \caption{$C_1$(left), \quad $C_5$(right)}
 \label{fig:C15}
\end{figure}

\begin{figure}[hbtp]
 \begin{center} 
 \includegraphics[width=6cm]{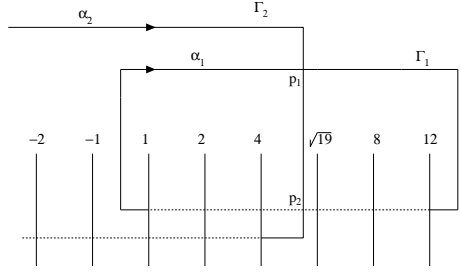}
  \end{center}
 \caption{$intersection$}
 \label{fig:int}
\end{figure}

 Next we construct another system that is dual to $L(\Gamma )$.
Let $r(i)$ be the oriented line segment from $s_0$ to $s_i$. Using this
notation we construct 8 elements
in $H_2(S_0,{\bf Z})$ ( see Fig.\ref{fig:C15}):
\begin{eqnarray*}
&&
C_1=(r(6)-r(7))\times \alpha _2,
\cr &&
C_2=(r(9)-r(10))\times \alpha _1,
\cr &&
C_3=(r(7)-r(12))\times \alpha _2,
\cr &&
C_4=(r(4)-r(9))\times \alpha _1,
\cr &&
C_5=r(10)\times \alpha _1-r(11)\times (\alpha _1+\alpha _2)
-r(12)\times \alpha _2,
\cr &&
C_6=r(4)\times \alpha _1-r(5)\times (\alpha _1+\alpha _2)-
r(6)\times \alpha _2,
\cr &&
C_7=-r(7)\times \alpha _2-r(8)\times (\alpha _1-\alpha _2)-
r(9)\times \alpha _1,
\cr &&
C_8=-r(1)\times \alpha _2-r(2)\times (\alpha _1-\alpha _2)-
r(3)\times \alpha _1.
\end{eqnarray*}
We note that any element of $\{ C_1,\cdots , C_8\} $ ends with a vanishing
cycle at the terminal point of the base arc. So it is a 2-cycle on $S_0$.
Let $L(C)$ denote the sublattice $\oplus {\bf Z}C_i$.
\begin{prop}
We have the the matrices of 
the intersection numbers $\Gamma _i \cdot \Gamma _j$ and 
$\Gamma_i \cdot C_j$ as follows :
\begin{eqnarray*}
&&
(\Gamma _i\cdot \Gamma _j)=\begin{pmatrix}
0&2&0&0&0&0&0&0\\
2&0&0&0&0&0&0&0\\
0&0&0&2&0&0&0&0\\
0&0&2&0&0&0&0&0\\
0&0&0&0&-2&0&0&0\\
0&0&0&0&0&-2&0&0\\
0&0&0&0&0&0&-2&0\\
0&0&0&0&0&0&0&-2 \end{pmatrix}.
\cr &&
(\Gamma _i\cdot C_j)=\begin{pmatrix}
1&0&0&0&0&0&0&0\\
0&1&0&0&0&0&0&0\\
0&0&-1&0&0&0&0&0\\
0&0&0&1&0&0&0&0\\
0&-1&-1&0&1&0&0&0\\
1&0&0&0&0&1&0&0\\
1&1&-1&-1&0&0&1&0\\
0&0&0&0&0&0&0&1 \end{pmatrix}.
\end{eqnarray*} 
\end{prop}
[Proof]. We can calculate all the intersection numbers by a direct observation
of the cycles. For example we consider $\Gamma _1\cdot \Gamma _2$. They have
2 intersections on a vertical line between $4$ and $\sqrt{19}$ of $\Gamma _2$
saying $p_1,p_2$ (see Fig.\ref{fig:int}).
At the upper intersection $p_1$ the base arcs intersect with a negative
sign, and the fibre cycles $\alpha _1,\alpha _2$ intersect with a
positive sign. By considering the definition of the orientation 
on $\Gamma _i$ we get there intersection number $+1$. The situation is 
quite the same for the intersection in the lower half plane. 
So we have $\Gamma _1\cdot \Gamma _2=2$. We can proceed this type
of calculation until we get the full intersection matrices.
\par
\rightline{q.e.d.} 
\begin{rem}
1)\quad By using the conventional notations $U=\begin{pmatrix}0&1\\ 1&0\end{pmatrix}$ and 
$U(k)=\begin{pmatrix}0&k\\ k&0\end{pmatrix}$, we have
\[
(\Gamma _i\cdot \Gamma _j)=
U(2) \oplus U(2) \oplus \langle -2 \rangle ^4.
\]
2)\quad We note that $(\Gamma _i\cdot C_j)$ is a unimodular matrix.
So we can find a generator system 
$\{ \check{\Gamma }_1,\cdots ,\check{\Gamma }_8\} $ of $L(C)$ with
$\Gamma _i\cdot \check{\Gamma }_j=\delta _{ij}$.
\end{rem}
\begin{lem}
$L(\Gamma )$ and $L(C)$ are nondegenerate rank 8 lattices. 
$L(\Gamma )$ is a primitive
sublattice in $H_2(S_0,{\bf Z})\cong L$.
\end{lem}
Proof].\quad The first statement is a direct consequence of the fact that 
$(\Gamma _i\cdot C_j)$ is a unimodular matrix.
\par
Suppose $L(\Gamma )$ is not primitive. It means that there is an element
$\lambda  \in L-L(\Gamma )$ with $k\lambda  \in L(\Gamma )$ for
some integer $k$.Put
\[
k\lambda =\sum_{i=1}^8k_i\Gamma _i.
\]
Then we have 
\[
k(\lambda \cdot \check{\Gamma}_i)=(k\lambda )\cdot \check{\Gamma}_i=k_i.
\]
It indicates $k\mid k_i\ (i=1,\cdots ,8)$. So it is deduced $\lambda \in L(\Gamma )$.
This is a contradiction.
\par
\rightline{q.e.d}
\par
\subsection{The algebraic sublattice}
\quad For the moment we study the divisors on the reference surface $S_0$.
Let $P_{ij}^{\pm }$¡¡denote two intersections 
$H_i\cap H_j \ \ (1\leq i<j\leq 4)$, 
here we distinguish
them by the signature of the $s$-coordinate ( see Proposition 1.1). 
The surface $S_0$ has an exceptional 
curve corresponding to every intersection $P_{ij}^{\pm }\ \ (1\leq i<j\leq 4)$. 
We denote them by $E_{ij}^{\pm }$.
Let us consider the projection $\pi:S \rightarrow {\bf P}^1\times {\bf P}^1$.
Let $F_s$ denote the pull back $\pi ^{\ast }({\bf P}^1(s))\times \{t\} )$ of the
generic $s$-coordinate line, and let $F_t$ denote the pull back of
the generic $t$-coordinate line.
Set
\[
G_i=\frac{1}{2}(\pi ^{\ast }H_i-\sum_{j\not= i}(E_{ij}^{+}+E_{ij}^{-})),
\ (i=1,2,3,4).
\]
This is the reduced divisor coming from the 4 sections $H_i$ of the
elliptic fibration $(S_0,\pi_1 ,{\bf P}^1(s))$.By an easy observation we have
the following.
\begin{lem}
We have intersection numbers among $E_{ij}^{\pm}$, $G_i$, $F_s$, $F_t$ :
\begin{eqnarray*}
&&
G_i\cdot G_i=-2,\ \ F_s\cdot F_t=2,\ \  G_i\cdot F_s=G_i\cdot F_t=1,
\cr &&
 G_i\cdot E_{ij}^+=G_i\cdot E_{ij}^-=1,\ \ 
 E_{ij}^+\cdot E_{ij}^+=E_{ij}^-\cdot E_{ij}^-=-2,
\end{eqnarray*}
and all the other intersections are 0.
We let $L(div)$ denote the sublattice of $H_2(S_0,{\bf Z})=L$ generated 
by $E_{ij}^{\pm},G_i,F_s,F_t$.
\end{lem}
\begin{rem}
If we observe the construction of $\Gamma _i$, 
we can show that $L(div)\subset L(\Gamma )^{\perp }$. Where $\perp $ indicates
the orthogonal complement in the full lattice $H_2(S_0,{\bf Z}) \cong L$.
\end{rem}
\begin{prop}
We have
\begin{eqnarray*}
&&
L(\Gamma )^{\perp}=L(div)=\langle E_{ij}^{\pm},G_i,F_s,F_t \rangle 
\cr &&
\cong (D_4)^3\oplus \langle -2\rangle \oplus \langle 2\rangle 
\end{eqnarray*}
on the reference surface $S_0$.
\end{prop}
Proof].\quad Set the sublattices $L_i\ \ (i=1,2,3)$ as follows:
\begin{eqnarray*}
&&
L_1=\langle G_1,E_{14}^{+},E_{14}^{-},F_t-E_{23}^+\rangle ,
\cr &&
L_2=\langle G_2,E_{24}^{+},E_{24}^{-},F_t-E_{13}^+\rangle ,
\cr &&
L_3=\langle G_3,E_{34}^{+},E_{34}^{-},F_t-E_{12}^+\rangle .
\end{eqnarray*}
They are isometric with the lattice 
\[
D_4=\begin{pmatrix}-2&1&1&1\\
             1&-2&0&0\\
             1&0&-2&0\\
             1&0&0&-2 \end{pmatrix} ,
\] 
and they are perpendicular
each others. Set
\begin{eqnarray*}
&&
\Delta _1=F_s+F_t-E_{12}^+-E_{13}^+-E_{23}^+,
\cr &&
\Delta _2=G_1+G_2+G_3-G_4+F_s+3F_t-E_{12}^+-E_{13}^+-E_{23}^+,
\end{eqnarray*}
then we have 
\[
\Delta _1\cdot \Delta _1=-2,
\Delta _1\cdot \Delta _2=0,
\Delta _2\cdot \Delta _2=2.
\]
Now we put
\[
P =L_1\oplus L_2\oplus L_3\oplus \langle \Delta _1\rangle
\oplus \langle \Delta _2\rangle  .
\]
Then it holds 
\begin{eqnarray*}
&&
P\subset L(div)\subset L(\Gamma )^{\perp },
\cr &&
{\rm rank }\ P=14, {\rm discr}\ P=-2^8.
\end{eqnarray*}
On the other hand we have ${\rm discr}\ L(\Gamma )=2^8$.
By considering the fact that $L$ is unimodular and $L(\Gamma )$ is
a primitive sublattice, we get the conclusion that
$P=L(div)=L(\Gamma )^{\perp }$.
\par
\rightline{q.e.d.}
Now we extend Proposition 1.4 to the general situation. 
\begin{thm}
Let $S$ be a double $4H$ surface. If we have ${\rm rank \ Pic}(S)=14$,
then it holds
\begin{eqnarray*}
&&
{\rm Tra}(S)\cong U(2)\oplus U(2)\oplus (-2I_4), 
\cr &&
{\rm Pic}(S)\cong D_4\oplus D_4\oplus D_4\oplus \langle -2 \rangle 
\oplus \langle 2 \rangle.
\end{eqnarray*}
\end{thm}
[Proof].\quad We note that the family ${\cal F}=\{ S(x)\} $ is fibred over the
parameter space of  ${\rm M}(2,{\bf C})^4$ . 
This fibration is locally 
topological trivial on a Zariski open subset in 
$X^{\circ} \subset {\rm M}(2,{\bf C})^4$. 
Using this trivialization
we can proceed the same argument for a general member as for the specialized
element $S_0$.In fact always $L(div)$ and $P$ are defined as a sublattice
in ${\rm Pic}(S(x))$. We have the systems $L(\Gamma )$ and $L(C)$
also. So Proposition 1.4 is valid for every $S$ of $\mathcal{F}$.
If the condition is satisfied, then $L(\Gamma )$ cannot contain any
divisor class. Hence we obtain the required conclusion.
\par
\rightline{q.e.d}
\begin{rem}
\par
As we will see in Section 2, we have ${\rm rank \ Pic}(S)=14$ for a
general member of ${\cal F}$.
\end{rem}

\par
\par

\section{Family ${\cal F}$ as that of lattice $K3$ surfaces}
\par
\subsection{General properties}
\quad
In this section we show the converse of Theorem 1. 
Before starting this argument
we state some general properties of a $K3$ surface.
\begin{lem}
Let $S$ be a $K3$ surface that is given by a minimal nonsingular model.
\par
1)\ \ For an irreducible curve $C$ on $S$, $C^2 \geq 0$ or $C^2 = -2$.
\par
2) If $D \in {\rm Pic(S)}$ satisfies\ \ $D^2 \geq -2$, 
then either $D$ or $-D$ is an effective class ( that is 
given by an effective divisor).
\end{lem}
This is deduced from the Riemann-Roch theorem by a routine argument.
\par
\begin{df}
We call an element $D\in {\rm Pic}(S)$ is {\em nef} when it holds
$D\cdot C\geq 0$ for any effective class $C$. 
\end{df} 

\begin{prop}[Pjatecki\u{i}-\v{S}apiro and \v{S}afarevi\v{c}\ \cite{PS}]
\ \
Let $S$ be a $K3$ surface, then we have the following:
\par
1)\ \ Suppose $x \in {\rm Pic}(S)$ satisfies\ \ $x \ne 0, \ \ x^2 = 0$. 
Then there exists an isometry $\gamma $ of ${\rm Pic}(S)$ such that $\gamma(x)$ becomes 
to be effective and nef. \par
2)\ \ Suppose $x \in {\rm Pic}(S)$ is effective, nef and $x^2 = 0$ then
$x$ is a multiple class of 
 an elliptic curve (i.e. $x = m[E]$ for a certain $m \in {\bf N}$ and an
 elliptic curve $E$ on $S$).\par
3)\ \ The linear system of an elliptic curve $C$ on $X$ determines an elliptic 
fibration $S \rightarrow {\bf P}^1$.
\end{prop}
\subsection{A lattice $K3$ surface and its realization as a 
double $4H$ surface}
\quad
Let us consider the lattice
\begin{eqnarray}
&&
P = {D}_4^3 \oplus \langle-2\rangle \oplus \langle2\rangle ,
\end{eqnarray}
and we let $(\ast ,\ast )$ denote the bilinear form on an abstract lattice $P$.
\begin{df}
We say a $K3$ surfaces $S$ is $\mathrm{of\ type} \ P$ if \par
$\mathrm{(a)} \ S$ admits an embedding $P \hookrightarrow \mathrm{Pic}(S)$,
\\
and $S$ is $\mathrm{of\ exact\ type}\ P$ if \par
$\mathrm{(a')} \ S$ admits an isomorphism $P \cong \mathrm{Pic}(S)$.
\end{df}
By the argument in Section 1 we can find the generator system of $P$:
\begin{eqnarray}
&&
f_s,\quad f_t,\quad g_i \ \ (i = 1,\cdots,4),\quad 
e_{ij}^+,\ \ e_{ij}^- \ \ (i \ne j, 1 \leq i,j \leq 4) , 
\end{eqnarray}
with the following properties among them :
\par
(p1)
\begin{eqnarray*}
&&
(g_i,g_i) = -2, \quad (f_s,f_t) = 2, \quad (g_i,f_s) = (g_i,f_t) = 1 
\cr &&
(g_i,e_{ij}^+) = (g_i,e_{ij}^-) = 1, \quad
 (e_{ij}^+,e_{ij}^+) = (e_{ij}^-,e_{ij}^-) = -2,
\end{eqnarray*}
and all other intersections equal $0$, 
\par
(p2)
\begin{eqnarray*}&&
f_s + f_t = 2g_i + \sum_{j \ne i}(e_{ij}^+ + e_{ij}^-).
\end{eqnarray*}
We note that this system is determined uniquely up to 
isometries of $P$.
\begin{thm}
Let $S$ be a $K3$ surface of type $P$ with the property:
\par
{\rm (n)} \ $f_s, f_t$ are nef classes under the identification in $\mathrm{(a)}$.
\\
Then $S$ has a representation as an extended double $4H$ surface.
\par
When $S$ is of exact type $P$ with the above conditioin ${\rm (n)}$, 
$S$ has a representation as a double $4H$ surface. 
\end{thm}

[Proof].
Always we regard the element of $P$ under the identification in ${\rm (a)}$ or
${\rm (a')}$.
We show the first assertion.
\par
(Claim 1) \quad We can assume that every $e_{ij}^{\pm }$ is effective.
\par
Let $e$ be one of the twelve elments $e_{ij}^{\pm}$. Since 
$e\cdot e = -2$, either $e$ or $-e$ is effective. If $-e$ is effective then we
perform the reflection $\gamma$ on ${\rm Pic}(S)$ determined by $e$. And we can use 
\[ \gamma(f_s),\quad \gamma(f_t),\quad \gamma(g_i),\quad 
    \gamma(e_{ij}^{\pm}) \]
instead of the system $(2.2)$. 
According to the orthogonality of the system $\{ e_{ij}^{\pm}\} $ 
we can iterate this procedure until we get the required effective system.
Since we have $\gamma(f_s) = f_s, \gamma(f_t) = f_t$, the
nef property for $f_s, f_t$ is always satisfied throughout this process. 
\par
\indent 
(Claim 2)\quad We can find a double covering $S'$, that is birationally equivalent with $S$,  
\[
\pi =(\pi _1,\pi _2): S'\rightarrow {\bf P}^1\times {\bf P}^1,
\]
ramified along a bidegree $(4,4)$ curve $B$. $S'$ has at most simple
singularities coming from the singular points of $B$. 
\par
By Proposition 2.1
$f_s, f_t$ are multiple classes of elliptic curve. Observing the assumption
$g_i\cdot f_s = g_i\cdot f_t = 1$ we know that they are reduced classes.
Again by Proposition 2.1 $f_s$ and $f_t$ determine 
two different elliptic fibrations $\pi_1 : X \rightarrow {\bf P}^1$ and
$\pi_2 : X \rightarrow {\bf P}^1$, respectively. Set
\[ 
\pi = (\pi_1,\pi_2) : S \rightarrow {\bf P}^1 \times {\bf P}^1 .
\]
This map is surjective and of degree 2, because we have
$\pi_1^{-1}(x)\cdot \pi_1^{-1}(y) = f_s\cdot f_t = 2$
for any $(x,y) \in {\bf P}^1 \times {\bf P}^1$. 
Let  $L_s, L_t$ be two lines ${\bf P}^1 \times \{\infty\}$ and 
$\{\infty\} \times {\bf P}^1$, respectively. 
So we obtain a double covering 
$\pi = (\pi_1,\pi_2) : S \rightarrow {\bf P}^1 \times {\bf P}^1 $. 
Generic fibers $\pi_1 ^{\ast }(x)$ and $\pi_2 ^{\ast }(y)$
are elliptic curves realized as double coverings over ${\bf P}^1$. 
So each of them has four branch points. Consequently, the branch locus $B$ of
$\pi$ is a curve of degree $(4,4)$. Since $S$ is a $K3$ surface, the canonical class
 $K_S = 0$. Hence $B$ has at most simple singularities.
\par
\indent (Claim 3)\quad  Let $G_i$ be the effective
divisor representing $g_i$.Then $H_i = \pi_*G_i$ is a $(1,1)$-curve $(i=1,2,3,4)$.
\par
Since $g_i\cdot g_i = -2$, either $g_i$ or $-g_i$ is effective. But $f_s$ is nef
 and has intersection $f_s\cdot g_i = 1$, then we know that $g_i$ is effective.  
 By the projection formula, we obtain
\[ 
\pi_*G_i\cdot L_s = G_i\cdot \pi^*L_s = g_i\cdot f_s = 1 .
\]
By the same argument, we get $\pi_*G_i\cdot L_t = 1$ also. 
Hence we obtain the required property.
\\
\indent (Claim 4)\quad
The effective class $e_{ij}^{\pm }$ is obtained by an exceptional divisor 
coming from the singularity of $B$.
\par
Let $E_{ij}^{\pm}$ be the effective divisor representing $e_{ij}^{\pm}$. 
Since $(e_{ij}^{\pm},f_s) = (e_{ij}^{\pm},f_t) = 0$
we get $\pi_*E_{ij}^{\pm} = 0$. 
It indicates that every  
$E_{ij}^{\pm}$ is an exceptional divisor derived from the singularity of $B$.
\\
\indent (Claim 5).\quad We have $B=H_1+H_2+H_3+H_4$
\par 
We consider the $(1,1)$ curve $H_i = \pi_*G_i$. 
As a divisor class $\pi^*H_i$ equals to $f_s + f_t$. 
By the starting assumption ${\rm (p2)}$, $\pi^*H_i$ and 
$2G_i + \sum_{j \ne i}(E_{ij}^+ + E_{ij}^-)$ are linearly equivalent. 
So there exist a principal divisor $D$ such that
\[ 
\pi^*H_i = 2G_i + \sum_{j \ne i}(E_{ij}^+ + E_{ij}^-) + D .
\]
Then 
\[ 
2H_i = \pi_*\pi^*H_i = 2\pi_*G_i + \pi_*D = 2H_i + \pi_*D .
\]
This is an equality as divisors itsselves (not as divisor classes), so we get 
$\pi_*D = 0$. This implies that $D$ is an sum of exceptional divisors. 
But $D$ is principal, so we must have $D = 0$. By observing the equality
\[ 
\pi^*H_i = 2G_i + \sum_{j \ne i}(E_{ij}^+ + E_{ij}^-) 
\]
we know that $\pi^*H_i$ has six components of the form $E_{ij}^{\pm}$,
note that we don't have any cancellation with an effective divisor $G_i$.  
Hence $H_i$ meets six different singular points on $B$. 
The sum of intersection numbers at these 
double points exceeds the number $(B,H_i) = 8$. 
In case $H_i$ is irreducible it must be an irreducible component of $B$.
In case we have $H_i=\ell _1+\ell _2$ with a $(0,1)$ curve and a $(1,0)$ curve,
at least one component, saying $\ell _1$, is contained in $B$. If $\ell _2$ is 
not contained in $B$, we have $\ell _2\cdot (B-\ell _1)=3$. 
So $\ell _2$ contains only one $\pi(E_{ij}^{\pm} )$.
 Namely $\ell _1$ contains 5 others.
Observing $\ell _1\cdot (B-\ell_1)=4$ we know that is too many. 
So we have $H_i\subset B$. Hence we obtained the claim.
\par
According to the above arguments we obtained the double covering
$\pi : S\rightarrow {\bf P}^1\times {\bf P}^1$. 
realizing the extended double $4H$ surface.
\par
Next let us prove the second assertion using the condition ${\rm (a')}$. 
It is enough to show  
\par
(Claim 6)\quad $S'$ has exactly 12 singular points of type $A_1$ on $\pi ^{\ast }B$.
\par
Let $L(exc)$ be the sublattice generated by the system $\{ e_{ij}^{\pm }\} $
in $P={\rm Pic}(S)$, this is isometric with 
$\langle -2\rangle ^{12}$. 
Suppose $x \in P$ is a class of irreducible
 $(-2)$-curve with $\pi_*x =0$. 
We can describe $x$ in the form 
\[
x = m_sf_s + m_tf_t + e + \sum_{i =1}^{4}m_ig_i , \quad
   e \in L(exc) 
\]
with the coefficients $m_s,m_t,m_i \in {\bf Z}$. 
By assumption it holds $\pi_*x\cdot L_s = \pi_*x\cdot L_t = 0$,
so we obtain
\[ 
m_s = m_t, \quad  2m_s + m_1 + m_2 + m_3 + m_4 = 0 .
\]
According to the starting condition ${\rm (p1), (p2)}$, we have
\begin{eqnarray*}
2(x - e) &=& 2(m_s(f_s + f_t) + m_1g_i + m_2g_2 + m_3g_3 + m_4g_4) \\
         &=& \sum_{i=1}^{4} m_i(2g_i - f_s - f_t) \\
         &=& -\sum_{i=1}^{4} m_i \sum_{j \ne i}(e_{ij}^+ + e_{ij}^-) 
\end{eqnarray*}
So we have $2x \in L(excep)$. 
By assumption we have $(2x,2x) = -8$. But such an element
 in $L(exc)$ should be the form $\pm2e_{ij}^{\pm}$ or 
$\pm e_1 \pm e_2 \pm e_3\pm e_4$. Where $\{ e_i\ \ (i=1,2,3,4)\} $ 
are distinct four
elements of  $\{ e_{ij}^{\pm}\} $.
In the case $2x = \pm2e_{ij}^{\pm}$, we conclude $x \in L(excep)$. 
\par
In the later case, we obtain 
\[ 
\pm e_1 \pm e_2 \pm e_3 \pm e_4- 2e = - 
    \sum_{i=1}^{4} m_i \sum_{j \ne i}(e_{ij}^+ + e_{ij}^-).
\]
Observe the right hand side, the number of the odd coefficients of $e_{ij}^{\pm}$ 
should be one of $0,6,8$. It does not attained in the left hand. So this
case does not happen.
Hence we obtained the claim. \par
\rightline{q.e.d}

\section{Modular group and marking with a multi-polarization}
\subsection{Congruence subgroup $G(2)$}
\quad
Let $L$ be the K3 lattice with the fixed basis as defined in Section 1 and $(\ ,\ )$ be 
the corresponding bilinear form. 
We consider the sublattice 
$P = (\mathrm{D}_4)^3 \oplus \langle -2 \rangle \oplus \langle 2 \rangle$
with the fixed generator system $\{ f_s, f_t, e_{ij}^{\pm }, g_i\} $ as in Section 2.
Set $T = P^{\perp} = \mathrm{U}(2)^2 \oplus \langle -2 \rangle^4$ with the fixed basis so that
$A=U(2)\oplus U(2)\oplus (-2I_4)$ be the intersection form. 
We denote the group of isometries of a lattice $M$ by ${\mathrm O}(M)$. Set $G = {\mathrm O}(T)$ and
\[ {\mathrm O}(L,P) = \{ g \in {\mathrm O}(L) : \quad g(x) = x \ \text{for} \ \forall
  x \in P \}.\]
Note that we have $g(T) = T$ for $g \in {\mathrm O}(L,P)$. So we can regard 
${\mathrm O}(L,P)$ as a subgroup of $G$. Now G is given as 
\[ G = \{ g \in M(8,{\bf Z}): \quad {}^tgAg = A \} .\]  
Set 
\[ G(2) = \{ g \in G: \quad g \equiv I \quad \text{mod} \quad 2 \}. \] 
Let $\alpha \in T$ be a $(-2)-$element, that is $(\alpha,\alpha) = -2$. 
It determines a reflection $\gamma _{\alpha }\in {\mathrm O}(L,P)$ by putting
\[ \gamma_{\alpha}(x) = x + (x,\alpha) \alpha .\] 
\begin{rem}
According to \cite{KKS} and \cite{MY}, we know that $G(2)$ is generated by
the reflections of the above form $\gamma_{\alpha }$. So it holds $G(2) \subset
{\mathrm O}(L,P)$.
\end{rem}
\begin{prop}\ We have
\[ G(2) = {\mathrm O}(L,P)\] . 
\end{prop}
[Proof]. \ 
It is enough to show ${\mathrm O}(L,P)\subset G(2)$. Here we use the notations and results in \cite{Ni}. 
Let $q_T:A_T\rightarrow {\bf Q}$ (resp. $q_P$) be the discriminant form of $T$ (resp. $P$),
where $A_T=T^{\ast }/T$.
\par
(Fact): $H=L/T\oplus P$ is a maximal isotropy subgroup for $q_T\oplus q_P$ in $A_T\oplus A_P$.
Namely it holds $(q_T\oplus q_P)\mid _H\equiv 0$ and $H$ is maximal with this property.
Moreover the canonical projections $H\rightarrow A_T$ and $H\rightarrow A_P$ induce the isomorphism
\[
A_T\cong H\cong A_P.
\]
\par
Hence we have 
\[
{\mathrm O}(L,P)\hookrightarrow  Ker\{ G\rightarrow {\mathrm O}(q_T)\}
\]
, where $G\rightarrow {\mathrm O}(q_T)$ indicates the canonical map. We
 can easily show that
\[
G(2)= Ker\{ G\rightarrow {\mathrm O}(q_T)\}.
\]
\par
\rightline{q.e.d}
\subsection{Marking with a multi-polarization}
\begin{df}
Let $S$ be a $K3$ surface with an embedding $P\hookrightarrow {\rm Pic}(S)$.
The triple $(S,\varphi, P)$ is a $P$-marking of $S$ provided 
\par
1) \ $\varphi : H_2(S,{\bf Z})\rightarrow L$ is an isometry,
 \par
2) \quad $\varphi ^{-1}(f_s),\varphi ^{-1}(f_t),\varphi ^{-1}(e_{ij}^{\pm}),\varphi ^{-1}(g_i)$ 
are effective and $\varphi ^{-1}(f_s),\varphi ^{-1}(f_t)$ are nef. 
\end{df}
According to Theorem 2 we have the double covering representation 
$\pi : S\rightarrow {\bf P}^1 \times {\bf P}^1$ and the class
\[ \pi_*(\varphi^{-1}(g_1) + \cdots + \varphi^{-1}(g_4)) \]
gives the class of the ramification divisor.
\begin{df}
Let $(S,\varphi,P)$ and $(S',\varphi',P)$ be $P$-markings of $S$ and $S'$, respectively.
An isomorphism $\rho : S \rightarrow S'$ is an isomorphism between these markings 
provided $\varphi = \varphi' \circ \rho_*$.
\end{df}
\begin{rem}
Because we have $\rho _{\ast }\varphi^{-1}(f_s)=(\varphi ')^{-1}(f_s)$ 
and $\rho _{\ast }\varphi^{-1}(f_t)=(\varphi ')^{-1}(f_t)$ ,
such an isomorphism $\rho$ preserves the covering structure over 
${\bf P}^1 \times {\bf P}^1$.
\end{rem}
 Namely we have the following commutative diagram 
with an element $\sigma \in 
\mathrm {PGL}_2(\bf C) \times \mathrm {PGL}_2(\bf C)$ :
\begin{eqnarray*}
\begin{CD}
  S      @>\text{$\rho$}>>    S'    \\
     @V\text{$\pi$}VV             @VV\text{$\pi'$}V   \\
 {\bf P}^1 \times {\bf P}^1  @>\text{$\sigma$}>>  {\bf P}^1 \times {\bf P}^1  
\end{CD}
\end{eqnarray*} 
\begin{df}
Let $H = (H_1,H_2,H_3,H_4)$, $H' = (H_1',H_2',H_3',H_4')$ be ordered
 sets of four curves of bidegree $(1,1)$. We say that $H$ and $H'$ are
 equivalent if there exists $\sigma \in 
\mathrm {PGL}_2(\bf C) \times \mathrm {PGL}_2(\bf C)$ such that $\sigma(H_i) = H_i' 
 \quad (i = 1,2,3,4)$.
\end{df}
\begin{rem}\ Let $(S,\varphi,P)$ and $(S',\varphi',P)$ be $P$-markings of $S$ and $S'$, respectively.
\par
(1)\ If these two $P$-markings are isomorphic then we have the same equivalent class of the ordered
sets $H=(H_1,H_2,H_3,H_4)$.
\par
(2)\ But we don't have the converse of (1). 
\end{rem}
\subsection{Modular group}
\quad 
Let $\{\Gamma_1, \cdots , \Gamma_8\}$ be the basis of $T$ such that 
$(\Gamma_i,\Gamma_j) = A$. We have elements 
$\{{\check \Gamma_1}, \cdots ,{\check \Gamma_8}\}$ in $L$ such that 
$({\check \Gamma_i},\Gamma_j) = \delta_{ij}$. 
The system  $\{ {\check \Gamma_i}\} $ is uniquely
determined modulo $P$. \par
Let $(S,\varphi,P)$ be a $P$-marking of a surface $S$. 
We consider the period
\[ [\int_{\varphi^{-1}({\check \Gamma_1})} \Omega: \cdots :
    \int_{\varphi^{-1}({\check \Gamma_8})}\Omega] \in {\bf P}^7 \]
where $\Omega$ is a holomorphic 2-form on $S$. The bilinear relation
\[ \int_S \Omega \wedge \Omega = 0, \quad 
   \int_S {\bar \Omega} \wedge \Omega > 0\]
implies that the period belongs to the domain
\[ D = \{\eta \in {\bf P}^7 : \quad {}^t \eta A \eta = 0, \quad 
   {}^t {\bar \eta}A\eta  > 0\}. \]
The domain $D$ has two connected components
\[ D = D^+ \cup D^-, \quad D^{\pm} = \{ (\eta_1, \cdots, \eta_8) \in D :
   \quad \pm{\rm Im}(\eta_3/ \eta_1) > 0 \}\]
and we can take $D^+$ as the period domain for the family of isomorphism classes of marked surfaces
$\{ (S(x), \varphi ,P):\ x\in X'\}$ .
\begin{rem}
Two domains $D^{\pm}$ are complex conjugate. That is
\[ [\int_{\varphi^{-1}({\check \Gamma_1})} \Omega: \cdots :
    \int_{\varphi^{-1}({\check \Gamma_8})}\Omega] \in D^+ 
   \quad \Leftrightarrow \quad 
   [\int_{\varphi^{-1}({\check \Gamma_1})} {\bar \Omega}: \cdots :
    \int_{\varphi^{-1}({\check \Gamma_8})} {\bar \Omega}] \in D^-. \]
\end{rem}
The group $G$ acts on the domain $D$. Set
\[ G(2)^+ = \{ g \in G(2) : \quad g(D^+) = D^+ \}.\]
It is a subgroup of $G(2)$ with index 2.
\begin{thm}
Let $(S,\varphi,P)$ and 
$(S',\varphi',l_s',P)$ be $P$-markings of extended double $4H$ surfaces $S$ and $S'$ , respectively.
Let $\eta, \eta' \in D^+$ be corresponding periods. Then these markings are isomorphic
 if and only if 
\[ g(\eta) = \eta' \]
for some $g \in G(2)^+\ (={\mathrm O}^+(L,P):=\{ g\in {\mathrm O}(L,P):\ g(D^+)=D^+\} )$.
\end{thm}
\begin{proof}
Assume $g(\eta) = \eta'$ for some $g \in G(2)^+$. 
According to Proposition 3.1 there exist an element ${\hat g} \in {\mathrm O}(L,P)$ such that 
${\hat g}|_{T} = g$. Then we have
\[ [\int_{\varphi^{-1} \circ {\hat g}({\check \Gamma_1})} \Omega: \cdots :
    \int_{\varphi^{-1} \circ {\hat g}({\check \Gamma_8})}\Omega] = 
   [\int_{(\varphi')^{-1}({\check \Gamma_1})} \Omega': \cdots :
    \int_{(\varphi')^{-1}({\check \Gamma_8})} \Omega'], \]
where $\Omega '$ indicates the holomorphic form on $S'$.
Consider the composition 
\[ f = (\varphi')^{-1} \circ {\hat g}^{-1} \circ \varphi : 
   {\rm H}_2(S,{\bf Z}) \longrightarrow L \longrightarrow L 
   \longrightarrow {\rm H}_2(S',{\bf Z}). \]
The above composite isomorphism induces the dual map 
\[ f^* : {\rm H}^2(S',{\bf Z}) \longrightarrow {\rm H}^2(S,{\bf Z})\]
with $f^*({\rm H}^{2,0}(S')) = {\rm H}^{2,0}(S)$. Moreover, $f^*$ 
preserves ample classes. Hence the Torelli theorem for K3 surfaces assures 
that there exists the unique isomorphism $\rho : S \rightarrow S'$ 
such that $\rho_* = f$. It is obvious that $\rho$ is an isomorphism of
 marked surfaces.\\
The converse is derived by the same argument.
\end{proof}

\setcounter{section}{3}
\section{The hypergeometric differential equation for the periods}
\quad
Let us take an element $C \in H_2(S_0,{\bf Z})$, and 
let $C(x) \in H_2(S(x),{\bf Z})$ denote its
continuation to $x=(x^1,x^2,x^3,x^4)\in X^0$ 
( that is multivalued and depends on the paths to $x$ in $X^{\circ}$). 
Now we investigate the differential
equation for the period
\begin{eqnarray}
u(x) =  \int_{C(x)} \Omega = 
\int \! \! \int _{C(x)} \{\prod _{p=1}^4(x_{11}^p\xi _1\eta _1+
x_{12}^p\xi _1\eta _2+x_{21}^p\xi _2\eta _1+x_{22}^p\xi _2\eta _2)^{-1/2}\}
\omega  \cr 
\omega = (\xi _1d\xi _2-\xi _2d\xi _1)\wedge (\eta _1d\eta _2-\eta _2d\eta _1),
\end{eqnarray}
where $[\xi _1,\xi _2]$ and $[\eta _1,\eta _2]$ denote the
homogeneous coordinates of ${\bf P}^1(s)$ and  ${\bf P}^1(t)$, respectively.
Note that this integral does not depend on the affine 
representatives of the homogeneous
coordinate of ${\bf P}^1(s)$ and  ${\bf P}^1(t)$. 
\par
There are left and right actions of $\mathrm {GL}(2,{\bf C})$ 
and a multiplicative action of $({\bf C}^{\ast })^4$ on 
$X^0$.
They induce the following actions on the period $u(x)$ :
\begin{eqnarray*}
&&
u(g\cdot x)=u(gx^1,\cdots ,gx^4) \quad g\in \mathrm {GL}(2,{\bf C}),
\cr \cr &&
u(x \cdot h)=u(x^1{}^th,\cdots ,x^4{}^th) \quad h\in \mathrm {GL}(2,{\bf C}),
\cr \cr &&
u(\lambda \circ x)=
u(\lambda _1x^1,\cdots ,\lambda _4x^4) \quad 
\lambda =(\lambda _1,\cdots ,\lambda _4)\in ({\bf C}^{\ast })^4
\end{eqnarray*}
\par \vskip 2mm
\begin{lem} We have the following equalities.
\begin{eqnarray}
&{\rm (1)}&\quad u(\lambda \circ x)=\prod _{p=1}^4\lambda _p^{a_p-1} u(x).
\\ \cr
&{\rm (2)}&\quad u(g\cdot x)=\frac{1}{{\rm det}(g)}u(x),\quad 
u(x\cdot h)=\frac{1}{{\rm det}(h)}u(x) .
\end{eqnarray}
\end{lem}
\begin{prop}  \label{proposition:pde}
\par
The integral $u(x)$
satisfies the following systems :
\begin{eqnarray}
\sum_{1\leq j,k \leq 2}x_{jk}^p\frac{\partial u}{\partial x_{jk}^p}&=&
-\frac{1}{2}u(x)\quad (p=1,2,3,4),
\\ \cr
E_1:\quad \sum_{p=1}^4\sum_{j=1}^2x_{lj}^p\frac{\partial u}
{\partial x_{mj}^p}&=&-\delta_{lm}u(x)\quad (\ell ,m\in \{ 1,2\} ),
\\ \cr
\sum_{p=1}^4\sum_{j=1}^2x_{jl}^p\frac{\partial u}{\partial x_{jm}^p}&=&
-\delta_{lm}u(x)\quad (\ell ,m\in \{ 1,2\} ),
\end{eqnarray}

\begin{eqnarray}
E_2:\quad 
\frac{\partial ^2u}{\partial x_{ij}^{q}\partial x_{k\ell }^p}
&=& \frac{\partial ^2u}{\partial x_{k\ell }^{q}\partial x_{ij}^p}
\quad (i,j,k,\ell \in \{ 1,2\} , p,q\in \{ 1,2,3,4\} ),
\\
\frac{\partial ^2u}{\partial x_{11}^{q}\partial x_{22}^p}
&=&\frac{\partial ^2u}{\partial x_{21}^{q}\partial x_{12}^p}
\quad (p,q\in \{ 1,2,3,4\} ).
\end{eqnarray}
\end{prop}
\par \vskip 3mm
[Proof]\quad 
Differentiate the first equality (4.2) in  Lemma 4.1 with respect to $\lambda _i$ and put
$\lambda =(\lambda _1,\lambda _2,\lambda _3,\lambda _4)=(1,1,1,1)$.
Then we get the first equality (4.4).
\par
Differentiate the second equality (4.3) in Lemma 3.1 for the left action 
with respect to the $ij$-component $g_{ij}$ of $g$ and put $g=I_2$.
Then we get the second equality (4.5). We get the third equality by 
using the equality for the right action of $\mathrm {PGL}(2,{\bf C})$ with the 
same method.
The system $E_2$ is deduced from the direct computation and 
the commutativity of the 
partial differentiations and the integral.
\par
\rightline{q.e.d.}
\par \vskip 3mm
Now we show the system $E_1+E_2$ is a holonomic system on $X^0$.
We consider the variety $B $ in the cotangent bundle
$T^{\ast }(X^{\circ })$ defined by 
\begin{eqnarray}
\sum_{1\leq j,k \leq 2}x_{jk}^p \xi _{jk}^p &=& 0 \quad (p=1,2,3,4),
\label{eq:linear1}
\\
\sum_{p=1}^4\sum_{j=1}^2x_{\ell j}^p\xi_{mj}^p &=& 0 \quad 
(\ell ,m\in \{ 1,2\} ), \\
\sum_{p=1}^4\sum_{j=1}^2x_{j\ell }^p\xi_{jm}^p &=& 0
\quad (\ell ,m\in \{ 1,2\} ),
\\ \cr
\xi_{ij}^{q}\xi_{k\ell }^p - \xi_{k\ell }^{q}\xi_{ij}^p &=&0
\quad (i,j,k,\ell \in \{ 1,2\} , p,q\in \{ 1,2,3,4\} ), 
\\ \cr
\xi_{11}^{q}\xi_{22}^p - \xi_{21}^{q}\xi_{12}^p &=&0
\quad (p,q\in \{ 1,2,3,4\} ),
\end{eqnarray}
where $\xi _{jk}^p$ stands for $\frac{\partial u}{\partial x_{ij}^k}$.
We regard $B$ as a fiber space over $X^{\circ }$ in the product of $X^{\circ }$ 
and a space of symbols
$\xi =(\xi ^1,\xi ^2,\xi ^3,\xi ^4)$ for
\[
\xi ^p=\begin{pmatrix}\xi_{11}^p&\xi_{12}^p \cr \xi_{21}^p &\xi_{22}^p\end{pmatrix}.
\]
The characteristic variety itself does not necessarily coincide with $B$, 
but it contains $B$, which we will call the fake characteristic variety. 
We will prove that the fake characteristic variety consists of only
$\{ 0 \}$ (zero-section) on $X^0$.
It implies that any solution of $E_1 + E_2$ on a simply connected
domain in $X^0$ is holomorphic.

\par \vskip 2mm
Let us fix an arbitrary point
\begin{eqnarray*}
x = (x^1,x^2,x^3,x^4)\in X^0, \quad x^p = \begin{pmatrix}
         x^p_{11} & x^p_{12}  \cr
         x^p_{21} & x^p_{22}\end{pmatrix} \quad (p=1,2,3,4). 
\end{eqnarray*}
Put
\begin{eqnarray*}
H_p : \quad f_p = x^p_{11} \xi_1 \eta_1 + x^p_{12} \xi_1 \eta_2 +
          x^p_{21} \xi_2 \eta_1  + x^p_{22} \xi_2 \eta_2 = 0 ,
\end{eqnarray*} 
and set
\begin{eqnarray*}
&&
D(pq)={\rm Det}(x^p)\ {\rm Trace}((x^p)^{-1}x^q)
=x_{11}^px_{22}^q+x_{22}^px_{11}^q-x_{12}^px_{21}^q-x_{21}^px_{12}^q.
\end{eqnarray*}
We set
\[
M_{pqr}=\begin{pmatrix}x_{11}^p&x_{12}^p&x_{21}^p&x_{22}^p\cr
                 x_{11}^q&x_{12}^q&x_{21}^q&x_{22}^q\cr
                 x_{11}^r&x_{12}^r&x_{21}^r&x_{22}^r\end{pmatrix},
 \quad p,q,r\in \{ 1,2,3,4\},
\]
and let $(ijk)$ denote the $3\times 3$ minor determinant of $M_{pqr}$ 
induced from the $i,j,k$-th column vectors. Put
\[
D(pqr)=(234)(123)-(134)(124).
\]
\begin{lem} \quad We have the following :
\par
{\rm (1)}\ $D(pq) = D(qp), \quad D(pp) = 2{\rm det}(x^p)$
\par
{\rm (2)}\ $H_p$  is irreducible if and only if 
   \[ D(pp) \ne 0 \]
\indent {\rm (3)}\ Let $H_p$ and $H_q$ be different and both irreducible . 
Then they have
different intersection points if and only if
\[ D(pq)^2 - D(pp) D(qq) \ne 0 .\] 
\par
{\rm (4)}\ We have $H_p\cap H_q\cap H_r=\emptyset$ if and only if 
$D(pqr)\not= 0$.
\end{lem}
\par \vskip 2mm
[Proof].
\par \indent
The first two claims are obvious. So we consider the third statement.
\par
$H_p$ is expressed in the form
\[ \frac{\xi_2}{\xi_1} = - \frac{x^p_{11} \eta_1 + 
     x^p_{12} \eta_2}{x^p_{21} \eta_1 + x^p_{22} \eta_2} . \]
So we obtain the intersections $H_p\cap H_q$ from
\[ \frac{x^p_{11} \eta_1 + x^p_{12} \eta_2}{x^p_{21} \eta_1 + 
   x^p_{22} \eta_2} = \frac{x^q_{11} \eta_1 + x^q_{12} \eta_2}
  {x^q_{21} \eta_1 + x^q_{22} \eta_2} .\]
Then the intersection comes from the eigen vector of $(x_p)^{-1}x^q$.
Hence we get the required condition.
\par
Next we examine the last statement. 
Generally we have
\[
M_{pqr}\begin{pmatrix}(234)\cr -(134)\cr (124)\cr -(123)\end{pmatrix}=\begin{pmatrix}0\cr 0\cr 0\end{pmatrix} .
\]
Note that we have $H_p\cap H_q\cap H_r\not= \emptyset $ if and only if
there exists a solution $\zeta ={}^t(\zeta _1,\zeta _2,\zeta _3,\zeta _4)$ with
$\zeta _1\zeta _4-\zeta _2\zeta _3=0$ for
\[
M_{pqr}\begin{pmatrix}\zeta _1\cr \zeta _2\cr \zeta _3\cr \zeta _4\end{pmatrix} =
\begin{pmatrix}0 \cr 0 \cr 0 \end{pmatrix}
\]
If ${\rm rank} M_{pqr}<3$ we can find easily such a solution $\zeta $.
In case  ${\rm rank} M_{pqr}=3$ we have 1 dimensional solution space
for $M_{pqr}\zeta =0$. So ${}^t((234),-(134),(124),-(123))$ becomes
a required solution only when $D(pqr)=0$.
\par
\rightline{q.e.d}
\begin{thm}  \label{theorem:non-singular}
Let $x=(x^1,x^2,x^3,x^4)\in \mathrm{M}(2,{\bf C})^4$ be a point on $X^0$, namely
$x$ satisfies the conditions:
\par
{\rm (g1)}. Any $H_p$ is irreducible i.e.  $D(pp) \ne 0$, 
\par
{\rm (g2)}. $H_p$ and $H_q$ have 2 different intersections i.e.
$D(pq)^2 - D(pp)D(qq) \ne 0$ for any $i\not= j$
\par
{\rm (g3)}. $H_p\cap H_q\cap H_r=\emptyset $ 
i.e
$D(pqr) \ne 0$ for any $(p-q)(q-r)(r-p)\not= 0$ \par
\noindent
Then any local solution of the system $E_1+E_2$ 
around the point $x \in X^0$ is locally holomorphic.
\end{thm}

[Proof]. We show $B_x=\{ 0\} $ in several steps.
\par \indent
(Step1 ).
Let us consider the Segre embedding
\[ 
\psi : {\bf P}^1 \times {\bf P}^1 \longrightarrow {\bf P}^3,
\]
that is defined by 
\[
[\xi_1,\xi_2] \times [\eta_1,\eta_2] \mapsto  
  \begin{pmatrix}\xi_1 \eta_1 & \xi_1 \eta_2 \cr 
  \xi_2 \eta_1 & \xi_2 \eta_2\end{pmatrix} .
\]
Note that the point 
${\displaystyle \begin{pmatrix}\xi_{11} &\xi_{12}\cr \xi_{21}& \xi_{22}\end{pmatrix} \in {\bf P}^3}$
belongs to ${\rm Im}(\psi )$ if and only if $\xi_{11}\xi_{22} - \xi_{12}\xi_{21}=0$.
Let $\Delta $ denote the diagonal map
\[ 
\Delta : {\bf P}^3 \longrightarrow {\bf P}^3 \times {\bf P}^3 
  \times {\bf P}^3 \times {\bf P}^3, \qquad P \mapsto (P,P,P,P). 
\]
Set
\[
\xi ^p=\begin{pmatrix}\xi_{11}^p&\xi_{12}^p\cr \xi_{21}^p&\xi_{22}^p\end{pmatrix} ,
\]
and regard $\xi =(\xi ^1,\xi ^2,\xi ^3,\xi ^4)$ as a homogeneous coordinate
on 
${\bf P}^3 \times {\bf P}^3 \times {\bf P}^3 \times {\bf P}^3$.
Then the system $E_2$ in Proposition 4.1
determines exactly the image of
$\Delta \circ \psi $. 
\par
(Step 2).
Let $\xi =(\xi ^1,\xi ^2,\xi ^3,\xi ^4)$ be a point on $B_x$.
Then at least two of $\{ \xi ^p\} $ should be $O$.
Suppose in contrary three of them , saying $\xi ^1,\xi ^2,\xi ^3$, 
are not equal to $O$. According to the argument in Step 1 they have the 
pull backs
in ${\bf P}^1\times {\bf P}^1$, and these are the same point 
$P=([\xi _1,\xi _2],[\eta _1,\eta _2])$.
Replace $\xi^p$ by $\psi(P)$ in (4.9),
which yields $ P \in H_p$.
Hence we obtain $P\in H_1\cap H_2\cap H_3$. This contradicts 
${\rm (g3)}$.
\par 
(Step 3).
Let $\xi =(\xi ^1,\xi ^2,\xi ^3,\xi ^4)$ be a point on $B_x$.
If three of $\{ \xi ^p\} $ are $O$, then all of them equal $O$.
\par \indent
Suppose $\xi ^2=\xi ^3=\xi ^4=0$. Then the equation (4.10) reduces to
\[
x^1{}^t\xi ^1=O.
\]
Because $x^1$ is assumed to be invertible, we obtain $\xi ^1=O$.
\par
(Step 4).
We don't have the case
$\xi^1 \ne 0, \xi^2 \ne 0, \xi^3 = \xi^4 = 0$.
\par \indent
Suppose it happens, then we have $\xi ^2=c\xi ^1$ for some constant
$c\ (\not= 0)$.
The two equations (4.10) and (4.11) are expressed in the form
\[ 
x^1 \  {}^t \xi^1 + x^2 \  {}^t \xi^2 = 0, \quad  
   {}^t \xi^1 \  x^1 + {}^t \xi^2 \  x^2 = 0 .
\]
Because we supposed $x^2$ to be invertible, we have
\[ 
(x^2)^{-1} \  x^1 \  {}^t \xi^1 = {}^t \xi^1 
    \  x^1 \  (x^2)^{-1}. 
\]
So we get 4 linear equations for 4 unknowns 
$\xi _{11}^1,\xi _{12}^1,\xi _{21}^1,\xi _{22}^1$. Together with the 2 equations
coming from the first equation (4.9) for $p=1,2$ we obtain the system of 
linear equations
\[ 
M \begin{pmatrix} 
   \xi^1_{11} \cr \xi^1_{12} \cr \xi^1_{21} \cr \xi^1_{22} \end{pmatrix}
   = \begin{pmatrix} 0 \cr 0 \cr 0 \cr 0 \cr 0 \cr 0  \end{pmatrix}\qquad \qquad (\ast ). 
\]
By a direct calculation we have
\[
M=
\begin{pmatrix}
x^1_{12} x^2_{21} - x^1_{21} x^2_{12} &
x^1_{12} x^2_{22} - x^1_{22} x^2_{12} & 
x^1_{22} x^2_{21} - x^1_{21} x^2_{22} & 
0 \cr
x^1_{11} x^2_{12} - x^1_{12} x^2_{11} &
0 &
x^1_{11} x^2_{22} - x^1_{22} x^2_{11} &
x^1_{12} x^2_{22} - x^1_{22} x^2_{12} \cr
x^1_{21} x^2_{11} - x^1_{11} x^2_{21} &
x^1_{22} x^2_{11} - x^1_{11} x^2_{22} &
0 &
x^1_{22} x^2_{21} - x^1_{21} x^2_{22} \cr
0 &
x^1_{11} x^2_{12} - x^1_{12} x^2_{11} &
x^1_{21} x^2_{11} - x^1_{11} x^2_{21} &
x^1_{21} x^2_{12} - x^1_{12} x^2_{21}\cr
x^1_{11} & x^1_{12} & x^1_{21} & x^1_{22} \cr
x^2_{11} & x^2_{12} & x^2_{21} & x^2_{22} 
\end{pmatrix} .
\]
Let $(i,j,k,\ell )$ denote the $4\times 4$ minor of $M$ obtained by
taking $i,j,k,\ell $-th row vectors. By assumption the system ($\ast $) has a
nontrivial solution $\xi ^1$, so every $(i,j,k,\ell )$ should be $0$.
By a direct calculation we have the following:
\begin{eqnarray*}
&&
(1256) = (x^1_{12} x^2_{22} - x^1_{22}x^2_{12}) (D(12)^2 - D(11) D(22)), \\ &&
(1356) = (x^1_{22} x^2_{21} - x^1_{21}x^2_{22}) (D(12)^2 - D(11) D(22)), \\ &&
(1456) = (x^1_{21} x^2_{12} - x^1_{12}x^2_{21}) (D(12)^2 - D(11) D(22)), \\ &&
(2356) = (x^1_{11} x^2_{22} - x^1_{22}x^2_{11}) (D(12)^2 - D(11) D(22)), \\ &&
(2456) = (x^1_{12} x^2_{11} - x^1_{11}x^2_{12}) (D(12)^2 - D(11) D(22)) ,\\ &&
(3456) = (x^1_{11} x^2_{21} - x^1_{21}x^2_{11}) (D(12)^2 - D(11) D(22)) ,
\end{eqnarray*}
and other $4\times 4$ minors are $0$.
Here recall the 2nd assumption for $x$. So we have $(D(12)^2 - D(11) D(22))$
is not $0$. Hence we obtain
\[
x_{ij}^1x_{k\ell }^2=x_{ij}^2x_{k\ell }^1
\]
for any indices $i,j,k,\ell \in \{ 1,2\} $. It means $x^2=c x^1$
and $H_1=H_2$. This is a contradiction.
\par
By the combination of the above arguments we obtain $B_x=\{ 0\} $.\\
\rightline{q.e.d.}

The system of partial differential equations $E_1 + E_2$
is closely related to the GKZ hypergeometric system
introduced by Gel'fand, Kapranov and Zelevinski.
Let us explain the relation and evaluate the dimension of the solution
space of $E_1+E_2$.
\par
Define a set of operators $E_1(\ell,m)$ and $E'_1(\ell,m)$ by
\begin{eqnarray*}
E(\ell,m) &:& \sum_{p,j} x^p_{\ell j} \partial^p_{m j} + \delta_{\ell m}, \\
E'(\ell,m)&:& \sum_{p,j} x^p_{j \ell} \partial^p_{j m} + \delta_{\ell m}. 
\end{eqnarray*}
Let $D$ be the Weyl algebra
$${\bf C}\langle 
      x^1_{11}, x^1_{12}, x^1_{21}, x^1_{22}, x^2_{11}, \ldots, x^4_{22},
      \partial^1_{11}, \partial^1_{12}, \partial^1_{21}, 
      \partial^1_{22}, \partial^2_{11}, \ldots, \partial^4_{22} \rangle.
$$               
Consider the GKZ-hypergeometric ideal
$H_A(\beta)$ in $D$
associated to the matrix 
$$
A =
 \left( \begin{array}{cccccccccccccccc}
1&1&1&1&0&0&0&0&0&0&0&0&0&0&0&0\\
0&0&0&0&1&1&1&1&0&0&0&0&0&0&0&0\\
0&0&0&0&0&0&0&0&1&1&1&1&0&0&0&0\\
0&0&0&0&0&0&0&0&0&0&0&0&1&1&1&1\\
1&1&0&0&1&1&0&0&1&1&0&0&1&1&0&0\\
1&0&1&0&1&0&1&0&1&0&1&0&1&0&1&0
\end {array} \right)
$$
and
$ \beta = (-1/2,-1/2,-1/2,-1/2,-1,-1)$.
\begin{prop}
Our system of partial differential equations $E_1$ and $E_2$
consists of three groups of operators
$$A \cdot \theta - {}^t \beta, $$
$$I_A=\{\partial^u-\partial^v \mid A u = A v, \ u, v \in {\bf N}_0^{16} \},$$
$$ E_1(\ell,m)\ {\rm and}\ E'_1(\ell,m)  (\ell \not= m). $$
Here, we denote by $\theta$ the column vector
$( x^i_{jk} \partial^i_{jk})$ of Euler operators of length 16.
\end{prop}
The GKZ-hypergeometric ideal $H_A(\beta)$ is the left ideal
in $D$ generated by
$A \cdot \theta- {}^t \beta$ and $I_A$.
Since the toric ideal $I_A$ is homogeneous, 
the $D$-module $D/H_A(\beta)$ is regular holonomic on $X={\bf C}^{16}$
(Hotta's theorem, see, e.g., \cite[p.82]{SST-BOOK}).
Our toric ideal $I_A$ satisfies the following properties:
\begin{enumerate}
\item The initial ideal of $I_A$ with respect to the reverse
lexicographic order is generated by square free monomials.
\item The toric ideal is Cohen-Macaulay.
\item The multiplicity of $I_A$ is $20$.
\item The variety $V(I_A)$ admits a natural action of $({\bf C}^*)^4$
and $V(I_A)/({\bf C}^*)^4$ is isomorphic to 
${\rm Im}\,(\Delta \circ \psi)$ in Theorem \ref{theorem:non-singular}.
\end{enumerate} 
The first statement can be easily checked by Buchberger's criterion.
The second fact follows from the first 
(see, e.g., \cite[p.153]{SST-BOOK}).
The third fact can be shown by computing the Hilbert polynomial of $I_A$
on computers.
Therefore, by theorems due to Gel'fand, Kapranov and Zelevinsky,
the rank of the solution space of $H_A(\beta)$
is $20$ and the singular locus agrees with the zero set of 
the principal $A$-determinant
(see, e.g., \cite[p.173]{SST-BOOK}).

We denote by E the left ideal in $D$ generated by
first order operators 
$E_1(\ell,m)$, $(\ell \not= m)$ and 
the GKZ-hypergeometric ideal $H_A(\beta)$.
\begin{thm}
The $D$-module $D/E$ is regular holonomic on $X$, and the rank of the solution space
of $E=E_1+E_2$ is equal to 8.
\end{thm}
[Proof]. 
Since $H_A(\beta) \subseteq E$ holds and the GKZ system is regular
holonomic, 
the $D$-module $D/E$ is also regular holonomic on $X$.
The differential operators $E_1(\ell,m), \ (\ell \not= m)$
are used to extract the space of the period maps from the $20$
dimensional solution space of the GKZ system.
\par
Next let us show the second statement.
Already in Section 3 we constructed 8 independent periods. So it is enough to say
that "at most 8 dimensional". 
Although it is possible in principle to evaluate the rank
by computer and Oaku's algorithm (see, e.g., \cite[p.31]{SST-BOOK}),
we could not evaluate it because of an exhaustion of memory.
So, we try to find sufficiently many initial terms
for the left ideal $E$ and a suitable weight.
For this system, we chose a weight
$ w^1 = 1 \cdot w, \, w^2 = 4 \cdot w, w^3 = 9 \cdot w, w^4 = 16 \cdot w$,
$ w = \begin{pmatrix}1 & 2 \cr
               0 & 4 \cr \end{pmatrix}$.
Put $W = (w^1, w^2, w^3, w^4)$.
The initial ideal generated by ${\rm in}_{(-W,W)}(E_1)$ and
${\rm in}_{(-W,W)}(E_2)$ has the rank $20$,
which means that we do not have sufficiently many initial terms.
We computed a partial Gr\"obner basis of $E_1$ and $E_2$ with the weight
$(-W,W)$ up to the degree 7 in the homogenized Weyl algebra.
The ideal generated by the initial terms of the partial Gr\"obner basis 
has rank $8$.
It follows from Theorem 2.5.1 of \cite{SST-BOOK} and the argument
of the regular holonomicity that the rank is bounded by $8$.
\par
\rightline{q.e.d.}
According to Theorem 4 and Theorem 5 the solution space of $E_1+E_2$ looks like
a vector bundle over $X^{\circ}$ of rank 8.
Let us take a $P$-marking $(S(x_0),\varphi ,P)$ , then we can choose 
the basis for the solution space of $E_1+E_2$
\[
\{ \int_{\varphi ^{-1}(\check{\Gamma}_1)}\Omega ,\cdots ,
\int_{\varphi ^{-1}(\check{\Gamma}_8)}\Omega \} 
\]
at $x_0\in X^{\circ}$.
So the system $E_1+E_2$ induces a representation of $\pi _1(X^{\circ },x_0)$
over $GL(8,{\bf Z})$.
\par
\begin{df}
The monodromy group {\rm Mono} for $({\cal F}, X^{\circ})$ is the
image of this representation.
\end{df}
\begin{rem}
Set $\Sigma =(GL(2,{\bf C})\setminus X^{\circ }/GL(2,{\bf C}))/({\bf C}^{\ast })^4$.
According to Lemma 4.1 we have the same period along the orbit
of the actions of $GL(2,{\bf C})$ and $({\bf C}^{\ast })^4$.
So the above monodromy representation reduces to that of $\pi _1(\Sigma ,\ast )$. 
\end{rem}
Let $B_r$ be the sublattice of $L$ generated by the elements
$f_s, f_t, g_i\ (i=1,\cdots ,4)$. Put
\[
O(L,Br)=\{ g\in O(L):\ g(P)=P, g\mid _{Br}={\rm id}\} .
\]
We define
\[
O^+(L,Br)=\{ g\in O(L,Br): g\mid _T(D^+)=D^+\} .
\]
The group $\pi _1(X^{\circ},x_0)$ acts on $L\cong H_2(S(x_0),{\bf Z})$ by the
natural way, and the branch locus is fixed under this action. So ${\rm Mono}$
is identified with a certain subgroup of $O^+(L,Br)$.

\begin{lem}
Let $(S_1,\varphi _1,P)$ and $(S_2,\varphi _2,P)$ be $P$-markings, and
let $x_1, x_2$ be corresponding points on $X^{\circ}$.
Suppose there is an isomorphism $\rho :S_1\mapsto S_2$ such that we have
$\varphi _2\circ \rho _{\ast }\circ \varphi _1^{-1}\in O^+(L,Br)$. 
Then the points $x_1$ and $x_2$ determine the
same point on $\Sigma =(GL(2,{\bf C})\setminus X^{\circ }/GL(2,{\bf C}))/({\bf C}^{\ast })^4$.
\end{lem}
[Proof].
\par
Because we have $\rho _{\ast}\varphi ^{-1}_1(f_s)=\varphi ^{-1}_2(f_s)$ and
$\rho _{\ast}\varphi ^{-1}_1(f_s)=\varphi ^{-1}_2(f_s)$, the isomorphism $\rho $
preserves the covering structure over ${\bf P}^1\times {\bf P}^1$. Moreover we have
$\rho _{\ast}\circ \varphi ^{-1}_1(g_i)=\varphi ^{-1}_2(g_i)\ (i=1,2,3,4)$. Hence
$\rho $ preserves the branch locus with the same numbering order. 
So $x_1$ and $x_2$ have the required property.
\par
\rightline{q.e.d.}
\par
\begin{thm}
We have
\[
O^+(L,Br)\mid _T={\rm Mono} .
\]
\end{thm}
\par
[Proof]. It is enough to show that $O^+(L,Br)\mid _T\subset {\rm Mono}$.
Let us fix a $P$-marking $(S_0,\varphi _0, P)$ corresponding to the
initial point $x_0\in X^{\circ}$.
Let $g$ be an arbitrary element of $O^+(L,Br)$, and put $g'=g\mid _T$.
Let $\eta _0\in D^+$ be the period point determined by $(S_0,\varphi _0, P)$,
and set $\eta _1= g'(\eta _0)$.
Let us take an oriented arc $r$ in $D^+$ that starts from $\eta _0$ and ends at $\eta _1$.
Because of the surjectivity of the period map, we can find a $P$-marking $(S(\eta ),\varphi (\eta ),P)$
for every point $\eta $ on $r$. Set $(S_1,\varphi _1, P)$ be the terminal marking and let
$x_1$ be the corresponding point on $X^{\circ }$.
Because of the injectivity of the period map we have
$(S_1,\varphi _1,P)=(S_1,\varphi _1\circ g^{-1},P)$.
Then by virtue of the above Lemma $x_0$ and $x_1$ determine the same point on $\Sigma $.
Hence we get unique closed arc $\gamma $ in $\Sigma $ corresponding to $r$.
It means that $g'$ is the monodromy transformation coming from the arc $\gamma $.
\par
\rightline{q.e.d.}
\setcounter{section}{4}
\section{Periods and Kuga-Satake varieties}
\quad 
In this section we construct the abelian variety attached to the
extended double $4H$ surface starting from its period.
The reader will find the method in \cite{vG} and \cite{S}. 
The detailed calculation and argument is exposed in \cite{Ko} also.
\par
Let us consider the lattice $T$ defined by the intersection matrix 
$A=U(2)\oplus U(2)\oplus (-2I_4)$ and
put $V_k=T\otimes k\  (k={\bf R}\ {\rm or}\ {\bf Q})$. Let $Q(x)$ denote the
quadratic form on $T$ and at the same time on $V_k$.
Let $Tens(T)$ and $Tens(V_k)$ be the
corresponding tensor algebras. And we let $Tens^+(T)$ and $Tens ^+(V_k)$
denote the subalgebras composed of the parts with even degree in
$Tens(T)$ and $Tens(V_k)$ respectively.
We consider the two sided ideal $I$ in 
$Tens ^+(V_k)$ generated by
the elements of the form $x\otimes x-Q(x)$ for $x\in V_k$, and the
ideal $I_{{\bf Z}}$ in $Tens(T)$ by the same manner. 
The corresponding even Clifford algebra is defined by
\[
C^+(V_k,Q)=Tens^+(V_k)/I.
\] 
By the same manner we define the even Clifford algebra over ${\bf Z}$ by
\[
C^+(T,Q)=Tens^+(T)/I_{{\bf Z}}.
\]
We note that $C^+(V_{\bf R},Q)$ is a $128$ dimensional real vector space and
$C^+(T,Q)$ is a lattice in it. So we obtain a real torus
\[
{\cal T}_{\bf R}=C^+(V_{\bf R},Q)/C^+(T,Q).
\]
Let ${\bf F}$ denote the quaternion algebra 
\[
{\bf Q}\oplus {\bf Q}i\oplus {\bf Q}j\oplus {\bf Q}ij
\]
with $i^2=j^2=-1$.
By some routine calculations of Clifford algebra we obtain the following.
\begin{prop}
We have the isomorphism of algebras 
$C^+(V_{\bf Q},Q)\cong  \mathrm{M}(4,{\bf F})\oplus \mathrm{M}(4,{\bf F})$.
\end{prop}
\par \vskip 2mm
Let a complex vector ${\underline \eta}=(\eta _1,\cdots ,\eta _8)$ 
be a representative of a point $\eta =[\eta _1 ,\cdots ,\eta _8]\in D^+$.
So it has an ambiguity of the multiplication by a non zero complex
number.
Put ${\underline \eta }=s+it\ \ (s,t\in {\bf R}^8)$.
If we impose the condition $(st)^2=-1$ in $C^+(V_{\bf R},Q)$, the representative
is uniquely determined up to a multiplication by a complex unit.
We denote it by
\[
{\underline \eta }=m_1(\eta )+im_2(\eta ),
\]
and put
\[
m(\eta )=m_1(\eta )m_2(\eta ).
\]
It is uniquely determined by $\eta $ without any ambiguity.
According to the imposed condition the element $m(\eta )\in C^+(V_{\bf R},Q)$ defines a complex structure on 
$C^+(V_{\bf R},Q)$
by the left action. 
It induces a complex structure on the real torus ${\cal T}_{\bf R}$ also.
We denote this complex torus by $(T,\ m(\eta ))$.
\par
Let $\{ \varepsilon _1,\cdots ,\varepsilon _8\} $ be the basis of $T$ with the 
intersection matrix $U(2)\oplus U(2)\oplus (-2I_4)$. And let 
$\{ e_1,\cdots ,e_8\} $ be a orthonormal basis of $V$ given by
\[
(e_1,\cdots ,e_8)=(\varepsilon _1,\cdots ,\varepsilon )
(\begin{pmatrix}\frac{1}{2}&\frac{1}{2}&0&0\cr
           0&0&\frac{1}{2}&-\frac{1}{2}\cr
           \frac{1}{2}&-\frac{1}{2}&0&0\cr
           0&0&\frac{1}{2}&-\frac{1}{2}\end{pmatrix} \oplus (I_4)) .
\]
Then the corresponding intersection matrix takes the form 
$I_2\oplus (-I_2)\oplus (-2I_4)$.
\par
Let $\iota $ be an involution on $C^+(V,Q)$ induced from the transformation
\[
\iota : e_{i_1}\otimes e_{i_2}\otimes \cdots \otimes e_{i_k}\mapsto 
e_{i_k}\otimes \cdots \otimes e_{i_2}\otimes e_{i_1}
\]
for the basis. Set $\alpha =4e_2e_1$.
According to the method in [St] we know that
\[
E(x,y)=tr(\alpha x^{\iota }y)
\]
determines a Riemann form on $( L,\ m(\eta ))$. We denote this abelian variety
by $A^+(\eta )$, that is so called the Kuga-Satake variety attached to the 
$K3$ surface corresponding to the period $\eta $.
In this way we can construct a family of abelian varieties 
\[
{\cal A}^+=\{ A^+(\eta )\ :\ \eta \in D^+\} 
\]
induced from
the lattice $T$ parameterized by the domain $D^+$. We can construct the
"conjugate family" 
\[
{\cal A}^-=\{ A^-(\eta )\ :\ \eta \in D^-\} 
\]
parameterized by 
\[
D^-=\{ \eta =[\eta_1,\cdots ,\eta_8]:\ {}^t\eta A\eta =0,\
{}^t{\overline \eta}A\eta >0,\ \Im (\eta_3/\eta_1)<0\}
\]
 by the same procedure with the Riemann form 
 $E^{-}(x,y)=-{\rm tr }(\alpha x^{\iota }y)$.
The right action of $C^+(V_{\bf Q},Q)$ on $(V,\ T,\ m(\eta ))$ commutes with the
left action of $\alpha (\eta )$. So we have
\[
C^+(T_{{\bf Q}}) \subset \rm{End}(A^{\pm }(\eta)) \otimes {\bf Q}  
\]
for any $A^{\pm }(\eta )$.
For a general member $\eta \in \cal{D}^+$, the endmorphism ring is given by 
\[ 
\rm{End}_{\bf Q}(A(\eta)) = \rm{End}(A(\eta)) \otimes {\bf Q} 
     \cong C^+(V_{{\bf Q}}). 
\]
According to Proposition 5.1 we obtain
\begin{thm}
For a general member $\eta \in \cal{D}^+$, $A^+(\eta)$ is isogeneous to a 
product of abelian varieties $(A_1(\eta) \times A_2(\eta))^4$ where 
$A_1(\eta)$ and $A_2(\eta)$ are 8-dimensional simple abelian varieties with 
$\rm{End}_{\bf Q}(A_i(\eta)) = {\bf F} \quad (i =1,2)$.
\end{thm}
\begin{rem}
Here we describe the relation between $A_1(\eta)$ and $A_2(\eta)$.
Now we define the linear involution $*$ on $V_{\bf R}$ by 
\[ e_1^* = - e_1 \quad {\rm and} \quad e_i^* = e_i 
        \quad (i = 2, \cdots, 8). \] 
It can be extended on $C^+(V_{\bf R},Q)$ as an autmorphism of algebra.
We define an involution $\sigma$ on $\cal{D}$ :
\[ \sigma : \mathcal{D} \longrightarrow \mathcal{D}, \quad
   (\eta_1, \cdots, \eta_8) \mapsto 
       (-\eta_2, -\eta_1, \eta_3, \cdots \eta_8). \]
So we have $\cal{D}_+^{\sigma} = \cal{D}_-$. It is easy to check that
we have
\[
A_2(\eta )\sim A_1(\eta ^{\sigma }),\ \ 
A_1(\eta )\sim A_2(\eta ^{\sigma }),
\]
where $\sim$ indicates the isogenous relation.
\end{rem}

\bigskip
\begin{flushright}
{\bf T. Tsutsui}

\medskip
Department of Mathematics and Informatics\\
Faculty of Science \\
Chiba University \\
Yayoi-cho Inage-ku Chiba 263-8522, JAPAN

\bigskip
{\bf K. Koike, H. Shiga}

\medskip
Graduate School of Science \\
Chiba University \\
Yayoi-cho Inage-ku Chiba 263-8522, JAPAN

\bigskip
{\bf N. Takayama}

\medskip
Department of Mathematics \\
Faculty of Science \\
Kobe University \\
Rokko Kobe 657-8501, JAPAN
\end{flushright}

\end{document}